\documentclass[a4paper,12pt,reqno]{amsart}

\usepackage[text={130mm,200mm}]{geometry}

\theoremstyle{plain}

\newtheorem*{theorem*}{Theorem}

\newtheorem*{corollary*}{Corollary}

\newtheorem*{lemma*}{Lemma}

\newtheorem*{proposition*}{Proposition}

\newtheorem*{conjecture*}{Conjecture}
\theoremstyle{definition}

\newtheorem*{definition*}{Definition}
\theoremstyle{remark}

\newtheorem*{remark*}{Remark}

\providecommand{\surname}[1]{#1}
\providecommand{\country}[1]{#1}
\begin{document}

\title[$G$-shtukas]{On  Langlands program,  related representation and $G$-shtukas}

\author{Nikolaj \surname{Glazunov}}
\address{ERMIT Department \\
  National Aviation University\\
  1~Komarova Pr.\\
  Kiev\\
 03680 \\
  \country{Ukraine}}
\email{glanm@yahoo.com}
\urladdr{https://sites.google.com/site/glazunovnm/}



% Submission date
%\date{Februaryr 29, 2020}
% \date{June 6, 2018; revised October 3, 2018}

\begin{abstract}
This communication is an introduction to the Langlands Program
and to ($G$-) shtukas (over algebraic curves) over function fields.
Modular curves and Drinfeld (elliptic) modules and shtukas are used in coding theory.
From this point of view the communication is concerned with mathematical models and methods 
of coding theory.
At first recall that linear algebraic groups found important applications in the Langlands program. Namely, for a connected reductive group $ G $ 
over a global field $ K $, the Langlands correspondence relates automorphic forms on $ G $ and global 
Langlands parameters,  i.e. conjugacy classes of homomorphisms  from the Galois group
 $ {\mathcal Gal} ({\overline K} / K) $ 
to the dual Langlands group $ \hat G ({\overline {\mathbb Q}} _ p) $. In the case of fields of 
algebraic numbers, the application and development of elements of the Langlands program made it 
possible to strengthen  the Wiles theorem on the Shimura-Taniyama-Weil hypothesis and to prove
 the Sato-Tate hypothesis.
 
  V. Drinfeld and V. Lafforgue have investigated the case of functional global fields of characteristic $ p> 0 $
   ( V. Drinfeld for $ G = GL_2 $ and V. Lafforgue for $ G = GL_r, \; r $ is an arbitrary positive integer). 
   They have proved in these cases the Langlands correspondence.
 
 Under the process of these investigations, V. Drinfeld  introduced the concept of a $ F $-bundle, or shtuka,
 which was used by both authors in the proof for functional global fields of characteristic $ p> 0 $ 
of the studied cases of the existence of the Langlands correspondence.
 
 Along with the use of shtukas developed and used by V. Drinfeld and V. Lafforge, other constructions
 related to approaches to the Langlands program in the functional case were introduced.
 
 G. Anderson has introduced the concept of a $ t $-motive. U. Hartl, his colleagues and students have
 introduced and have explored the concepts of finite, local and global $ G $ -shtukas.
 
 In this review article, we first present results on Langlands program and  related representation over algebraic number fields.
 Then we  briefly present approaches by U. Hartl, his colleagues and pupils to the study
 of $ G $ --shtukas.
 These approaches and our discussion relate to the Langlands program as well as to the internal development
  of the theory of $ G $-shtukas.

\end{abstract}

\maketitle

%%%%%%%%%%%%%%%%%%%%%%%%%%%%%%%%%%%%%%%%%%%%%%%%%%%%%%%%%%%%%%%%%%%%%%%%%%%%%%%%

\section*{Introduction}

This communication is an introduction to the Langlands Program
and to ($G$-) shtukas (over algebraic curves) over function fields.
Modular curves and Drinfeld (elliptic) modules and shtukas are used in coding theory.
From this point of view the communication is concerned with mathematical models and methods 
of coding theory.
We plan to consider applications of these models and methods to designing and researching 
of error correction codes  elsewhere. \\
At first recall that linear algebraic groups found important applications in the Langlands program.
 Namely, for a connected reductive group $ G $ 
over a global field $ K $, the Langlands correspondence relates automorphic forms on $ G $ and global 
Langlands parameters,  i.e. conjugacy classes of homomorphisms  from the Galois group
 $ {\mathcal Gal} ({\overline K} / K) $ 
to the dual Langlands group $ \hat G ({\overline {\mathbb Q}} _ p) $ \cite{Langlands1979,JacLan1970}. In the case of fields of 
algebraic numbers, the application and development of elements of the Langlands program made it 
possible to strengthen
 the Wiles theorem on the Shimura-Taniyama-Weil hypothesis and to prove the Sato-Tate hypothesis.
 Langlands reciprocity for $GL_n$ over non-archimedean local fields of characteristic zero is given by 
 Harris-Taylor \cite{HarrisTaylor}.
 
 V. Drinfeld \cite{Drinfeld1980} and V. Lafforgue \cite{Laforgue2002} have investigated the case of functional global fields of characteristic $ p> 0 $
   ( V. Drinfeld for $ G = GL_2 $ and V. Lafforgue for $ G = GL_r, \; r $ is an arbitrary positive integer). 
   They have proved in these cases the Langlands correspondence.
 
 In the process of these studies, V. Drinfeld introduced the concept of a $ F $ -bundle, or shtuka,
 which was used by both authors in the proof for functional global fields of characteristic $ p> 0 $ 
of the studied cases of the existence of the Langlands correspondence \cite{Drinfeld1987}.
 
 Along with the use of shtukas developed and used by V. Drinfeld and V. Lafforge, other constructions
 related to approaches to the Langlands program in the functional case were introduced.
 
 G. Anderson has introduced the concept of a $ t $-motive \cite{Anderson1986}. U. Hartl, his colleagues and postdoc students have
 introduced and have explored the concepts of finite, local and global $ G $ -shtukas\cite{Hartl2005,HartlViehmann2011,HartlRad2014,Singh2012,Rad2012,Weiss2017}.
 
 In this review, we first present results on Langlands program and  related representation over algebraic number fields.
 Then we  briefly present approaches by U. Hartl, his colleagues and pupils to the study
 of $ G $ --shtukas.
 These approaches and our discussion relate to the Langlands program as well as to the internal development
  of the theory of $ G $-shtukas.
  Some results on commutative formal groups and commutative formal schemes can be found in  
   \cite{Glazunov2015A,Glazunov2015B,Glazunov2015C} and in references therein.
  
 The content of the paper is as follows: \\
  1. Some results of the implementation of the Langlands program for fields of algebraic numbers and their localizations. \\
  2. Finite $ G $ -shtukas. \\
  3. Local $ G $ -shtukas. \\
  4. Global $ G $ -shtukas. \\

\section{Some results on Langlands program over algebraic number fields and their localizations} 

    Langlands conjectured that some symmetric power $ L$-functions extend to an entire function and coincide with certain automorphic $L$-functions.
    
\subsection{ Abelian extensions of number fields }

     In the case of algebraic number fields Langlands conjecture (Langlands correspondence) is the  global class field theory:    \\
          Representations of the abelian Galois group $Gal(K^{ab}/K)$ = characters of the Galois group $Gal(K^{ab}/K) $    
\begin{center}
                                                     correspond to  
\end{center}
automorphic forms  on $GL_1 $ that are characters of the class group of ideles. Galois group    $ Gal(K^{ab}/K)$ is the profinite completion of the group ${\mathbb A}^* (K)/K^* $ where $ {\mathbb A}(K) $ denotes the adele ring of $ K$.\\
  If $K$ is the local field, then Galois group $ Gal(K^{ab}/K)$ is canonically isomorphic to the profinite completion of $K^*$.   
    
\subsection{ $l$ - adic representations and Tate modules}

       Let $K$ be a field and $\overline K $  its separate closure, $ E_n=\{P \in E(\overline K )| nP=0\}$ the group of points of elliptic curve $E(\overline K)$  order dividing n. When $char K $ does not divide $n$ then $E_n$ is a free ${\mathbb Z}/n{\mathbb Z}$ -module of rank $2$.
       
             Let $l$ be prime, $l \ne char K$. The projective limit $T_l (E)$ of the projective system of modules $E_{l^m}$ is free
${\mathbb Z}_l$-adic Tate module of rank $2$.  

Let              
$V_l (E)=T_l (E) \otimes_{{\mathbb Z}_l} {\mathbb Q}_l$.
      Galois group $Gal({\overline K}/K)$ acts on all $E_{l^m}$ , so there is the natural continuous representation       ($l$-adic representation)
\begin{equation*}
\rho_{E,l} : Gal({\overline K}/K) \to Aut \; T_l (E) \subseteq Aut \;  V_l (E).  
\end{equation*}
$V_l (E)$ is the first homology group that is dual to the first cohomology group of $ l$-adic cohomology of elliptic curve $E$ and Frobenius $F$ acts on the homology and dually on cohomology. The characteristic polynomial $ P(T)$ of the Frobenius not depends on the prime number $ l$.

\subsection{Zeta functions and parabolic forms}

Let (in P. Deligne notations) $X$ be a scheme of finite type over ${\mathbb Z}$, $|X|$ the set of its closed points, and for each $x \in |X|$ let  $N(x)$ be the number of points of the residue field $k(x)$ of $X$ at $x$. The Hasse-Weil zeta-function of $X$ is, by definition
\begin{equation*}
\zeta_X (s)= \prod_{ x \in|X|}(1 - N(x)^{-s} )^{-1} .
\end{equation*}
In the case when   $X$  is defined over finite field ${\mathbb F}_q$, put  $ q_x= N(x)$, $deg(x)=[k(x):{\mathbb F}_q ]$, so $q_x=q^{deg(x)}$ . Put $t=q^{-s}$. Then 
\begin{equation*}
Z(X ,t) = \prod_{ x \in|X|}(1 - t^{deg(x)} )^{-1} .
\end{equation*}

The Hasse-Weil zeta function of $ E$ over ${\mathbb Q}$ (an extension of numerators of $\zeta_E (s)$ by points of bad reduction of $E$) is defined over all primes $p$:
\begin{equation*}
L(E({\mathbb Q}),s) = \prod_p (1 - a_pp^s + \epsilon(p)p^{1 - 2s} )^{-1} ,
\end{equation*}
here $\epsilon(p) = 1$ if $E$ has good reduction at $ p$, and $\epsilon(p) = 0$ otherwise.

For $GL_2 ({\mathbb R})$, let $C$ be its center,  $O(2)$ the orthogonal group.

Upper half complex plane has the representation: ${\mathbb H}^2=GL_2 ({\mathbb R}) / O(2) C$. So it is the homogeneous space of the group  $GL_2 ({\mathbb R})$.

A cusp (parabolic) form of weight $k \ge 1$ and level $N \ge1$ is a holomorphic function $f$ on the upper half complex plane ${\mathbb H}^2$ such that \\
a) For all matrices

$$g = \left(
\begin{array}{cc}
a& b \\
c& d
\end{array}
\right),     a,b,c,d \in {\mathbb Z}, a \equiv1(N),d\equiv1(N),c\equiv0(N) $$

and for all  $z \in {\mathbb H}^2$ we have
\begin{equation*}
f(gz)=f((az+b)/(cz+d))=(cz+d)^k f(z)
\end{equation*}
(automorphic condition).\\
b) 
\begin{equation*}
	|f(z)|^2 (Im z)^k 
\end{equation*}
is bounded on ${\mathbb H}^2$ . \\

 Mellin transform $ L(f,s)$ of the parabolic form $ f$ coincides with Artin $ L$-series of the representation $\rho_f$.
 
 The space ${\mathcal M}_n(N)$  of cusp forms of weight $ k$ and level  $N$ is a finite dimensional complex vector space. If 
 $f \in {\mathcal M}_n (N)$, then it has expansion 
 \begin{equation*}
f(z)=\sum_{n=1}^{\infty} c_n (f)\exp(2\pi inz)
\end{equation*}
and $L$-function is defined by
 \begin{equation*}
 L(f,s)=\sum_{n=1}^{\infty} c_n (f) /  n^s .
\end{equation*}

\subsection{Modularity results} 
  
  The compact Riemann surface $\Gamma\backslash{\mathbb H}^2$ is called the modular curve associated to the subgroup of finite index 
$\Gamma$ of  $GL_2 ({\mathbb Z})$ and is denoted by $X(\Gamma)$. If the modular curve is elliptic it is called the elliptic modular curve. \\    
The modularity theorem states that any elliptic curve over ${\mathbb Q}$ can be obtained via a rational map with integer coefficients from the elliptic  modular curve. \\
      By the Hasse-Weil conjecture (a cusp form of weight two and level $ N$ is an eigenform (an eigenfunction of all Hecke operators)). The conjecture follows from the modularity theorem.\\
         Recall the main (and more stronger than in Wiles \cite{Wiles} and in Wiles-Taylor \cite{TaylorWiles} papers) result by C. Breuil, B. Conrad, F. Diamond, R. Taylor \cite{BreuilConradDiamondTaylor}.
 
\begin{theorem*}
({\bf Taniyama-Shimura-Weil conjecture - Wiles Theorem.})
For every elliptic curve $E$ over ${\mathbb Q}$ there exists $f$, a cusp form of weight 2 for a subgroup 
$\Gamma_0 (N)$, 
such that $ L(f,s)=L(E({\mathbb Q}),s)$. 
\end{theorem*}
Recall that for projective closure $\overline E$ of the elliptic curve $E$ we have
\begin{equation*}
\overline E({\mathbb F}_p) = 1 - a_p +p. 
\end{equation*}
By H. Hasse
\begin{equation*}
a_p = 2\sqrt p \cos \varphi_p.
\end{equation*}
\begin{conjecture*}
({\bf Sato-Tate conjecture })
            How angles $\varphi_p$ are distributed in the interval $[0,\pi]$? Let $E$ be the elliptic curve has no complex multiplication. Then Sato have computed and Tate gave theoretical evidence that angles $\varphi_p$ in the case are equidistributed in $[0,\pi]$ with the Sato-Tate density measure  $\frac{2}{\pi} \sin^2 \varphi.$
\end{conjecture*}

The current state of Sato-Tate conjecture is now Clozel--Harris--Shepherd-Barron--Taylor Theorem 
\cite{Clozel HarrisTaylor,HarrisShepherd-BarronTaylor}.

\begin{theorem*}
(Clozel, Harris, Shepherd-Barron, Taylor). Suppose $E$ is an elliptic curve over ${\mathbb Q}$ with non-integral $j$ invariant. Then for all $n > 0; \; L(s;E; Sym^n)$ extends to a meromorphic function which is holomorphic and non-vanishing for $Re(s) \ge 1 + n/2$.
\end{theorem*}

These conditions and statements are sufficient to prove the Sato-Tate conjecture.

Under the prove of the Sato-Tate conjecture the Taniyama-Shimura-Weil conjecture oriented methods of A. Wiles and R. Taylor are used.

     Recall also that the proof of Langlands reciprocity for $GL_n$ over non-archimedean local fields of characteristic zero is given by Harris-Taylor \cite{HarrisTaylor}.
     
\section{Elliptic modules and Drinfeld shtukas.}   

%\newpage
\section{Finite $ G $ -shtukas.}
We follow to \cite{Drinfeld1987,Hartl2005,HartlViehmann2011,Singh2012}.
We start with very short indication on the general framework of the section. In connection with Drinfeld’s constructions of elliptic A-modules Anderson \cite{Anderson1986}  has introduced abelian t-modules and the dual notion of t-motives. Beside with mentioned papers these are the descent theory by A. Grothendieck \cite{Grothendieck}, cotangent complexes  by Illusie \cite{Illusie}, by S. Lichtenbaum and M. Schlessinger \cite{LichtenbaumSchlessinger}, by Messing \cite{Messing}  and by Abrashkin \cite{Abrashkin}. In this framework to any morphism $f: A \to B$ of commutative ring objects in a topos is associated a cotangent complex $ L_{(B/A)}$ and to any morphism of commutative ring objects in a topos of finite and locally free $Spec (A)$- group scheme $G$ is associated a cotangent complex $L_{(G/Spec (A))}$  as has presented in books by  Illusie \cite{Illusie}.
\subsection{Shtukas} 
Let $S$ be a scheme over $Spec \; {\mathbb F}_q$. 
In the paper \cite{Singh2012} author investigates relation between finite shtukas and strict finite flat commutative group schemes and relation between divisible local Anderson modules and formal Lie groups. Let $Nilp_{{\mathbb F}_{q}[[\xi]]}$ be the category of ${\mathbb F}_{q}[[\xi]]$-schemes on which $\xi$ is locally nilpotent. Let $S \in Nilp_{{\mathbb F}_{q}[[\xi]]}$. 
In Chapter 1 the author defines cotangent complexes as in papers by S. Lichtenbaum and M. Schlessinger  \cite{LichtenbaumSchlessinger}, by W. Messing \cite{Messing}, by V. Abrashkin \cite{Abrashkin} and prove that they are homotopically equivalent.\\
In section 1.5 the author  investigates the deformations of affine group schemes follow to the mentioned paper of Abrashkin and defines strict finite $\mathcal O$-module schemes.\\
Next section is devoted to relation between finite shtukas by V. Drinfeld \cite{Drinfeld1987} and strict finite flat commutative group schemes.\\
   The comparison between cotangent complex and Frobenius map of finite ${\mathbb F}_p$-shtukas is given in section 1.7.\\
     $ z$-divisible local Anderson modules by  Hartl \cite{Hartl2005}  and local schtukas are investigated in Chapter 2. \\
           Sections 2.1, 2.2 and 2.3 on formal Lie groups, local shtukas and divisible  local Anderson-modules define and illustrate notions for later use. Many of these, if not new, are set in a new form, 
       In Section 2.4 the equivalence between the category of effective local shtukas over $S$ and the category of $z$-divisible local Anderson modules over $S$ is treated. \\
         In the last section the theorem about canonical ${\mathbb F}_p[[\xi]]$ -isomorphism of $z$-adic Tate-module of $z$ -divisible local Anderson module $G$ of rank $r$ over $S$ and Tate module of local shtuka over $S$ associated to $G$ is given.
The main result of this dissertation is the following (section 2.5) interesting result
 it is possible to associate a formal Lie group to any $z$-divisible local Anderson module over $S$ in the case when $\xi$ is locally nilpotent on $S$.\\
Related and more general results have presented in the paper by U. Hartl, E. Viehmann \cite{HartlViehmann2011}.

\section{ Local $ G $ -shtukas}
Here we follow to \cite{HartlRad2014,Rad2012,HartlViehmann2011}.
Recall that local ${\mathbb P}$-shtukas are the functional field analogs of $p$-divisible groups with additional structure and moduli 
stacks of global $G$-shtukas are the functional field analogs for Shimura varieties.
Here ${\mathbb P}$ is a flat affine group scheme of finite type over $Spec \; {\mathbb F}[[z]]$ and $G$ is a flat affine group 
scheme of finite type over a smooth projective geometrically irreducible curve over ${\mathbb F}_{q}$; authors let
 ${\mathbb F}$ be a 
finite field extension of  a finite field ${\mathbb F}_{q}$ with $q$ elements and characteristic $p$.
 
The authors main theme is the relation between global $G$-shtukas and local ${\mathbb P}$-shtukas and authors of the paper 
under review prove the analog of a theorem of Serre and Tate over functional fields `stating the equivalence between the deformations of 
a global $G$-shtuka and its associated local $\mathbb{P}$-shtukas` .
Results explaining relation between local $\mathbb{P}$-shtukas and Galois representations as well as results on the existence of spaces by 
Rappoport-Zink \cite{RappoportZink} and their use 
to uniformize the moduli stack of  global $G$-shtukas also presented.
 
Similar partial results have appeared in dissertation by Singh \cite{Singh2012}  and in dissertation by Arasteh Rad \cite{Rad2012} both written under U. Hartl (M\"unster).

The generalization of purity of the Newton stratification in the context of local $G$-shtukas with applications to a generalization of 
Grothendieck’s conjecture on deformations of $p$-divisible groups with given Newton polygons were obtained by Veihmann 
\cite{Veihmann}.

The range of topics of the paper\cite{HartlRad2014} is indicated by the titles of the sections: 1. Introduction; 2. Local $\mathbb{P}$-shtukas and 
global $G$-shtukas; 
3. Tate modules for local $\mathbb{P}$-shtukas; 4. The Rapoport-Zink spaces for local $\mathbb{P}$-shtukas; 5. The relation between global 
$G$-shtukas and local $\mathbb{P}$-shtukas. 

Section 1 of \cite{HartlRad2014} comprises main results and list of notation and conventions.
In section 2, elements of the theory of loop groups, of local and global shtukas are introduced. These include formal torsor, local $\mathbb{P}$-shtukas, 
groupoid  of local $\mathbb{P}$-shtukas, quasi-isogeny between two local $\mathbb{P}$-shtukas, global $G$-shtuka. Main results include 
Proposition 2.4 on a canonical equivalence between the category fibered in groupoids that assigne to each $\mathbb{F}$-scheme $S$ the groupoid 
consisting of all formal torsors (the generalization of results by Hartl, Viehmann \cite{HartlViehmann2011}
and Proposition 2.11 on rigidity of quasi-isogenies for local $\mathbb{P}$-shtukas.
Section 3 concerns the relation between local $\mathbb{P}$-shtukas and Galois representations which is given by the associated Tate module. 
The (dual) Tate functor and the rational (dual) Tate functor are defined. Under certain conditions (dual) Tate functors are equivalences between 
category of {\`e}lale local shtukas over $S$ to the category of finite free $\mathbb{F}[[z]]$-modules equipped with continuous action of the 
algebraic fundamental group of $S$ at its geometric point. 
Section 4 works out unbounded Rapoport-Zink spaces for local $\mathbb{P}$-shtukas.  Boundedness conditions for local $\mathbb{P}$-shtukas 
and bounded local $\mathbb{P}$-shtukas also introduced and investigated. Representability of the bounded Rapoport-Zink functor is proved. 
The results are too technical to be stated here in details. They are applied later to  uniformize the moduli stack of  global $G$-shtukas. 
 The authors conclude this impressive work by proving the rigidity of quasi-isogenies for global $G$-shtukas and by proving 
the analog of the Serre-Tate theorem.

\section{ Global $ G $ -shtukas}
Here we follow to \cite{HartlRad2014,Rad2012,HartlViehmann2011,Weiss2017}.

Let $C$ be a smooth projective geometrically irroducible curve with the
function field  ${\mathbb F}_q(C)$
 over a finite field ${\mathbb F}_q$ with $q$ elements.
 Let $G$ be a parahoric Bruhat-Tits group scheme over $C$.
 
 Author of the paper \cite{Weiss2017}  considers  ``a foliation structure for Newton strata moduli spaces of bounded global $G$-shtukas with $H$-level structure for an arbitrary parahoric Bruhat-Tits group $G$'' and  ``Igusa varieties''. 
 She obtaines a morphism (Main Theorem 0.1) to the moduli space of  global $G$-shtukas.
 The author then relates here foliation structure to Oort`s foliations, to  Harris and Taylor   and to Mantovan.
These results, although difficult to explane in a short reviiew, are well summerised in a short Introduction.
Below bounded global $G$-shtukas with $H$-level structure are considered.
Briefly, the general idea is to start with a foliation stucture on the moduli space of such global $G$-shtukas and describe
it ``as a product of a covering of central leaves by Igusa varieties with truncated Rapoport-Zink spaces''.
The Main Theorem 0.1 gives the morphism from the product of author`s Igusa varieties and trancated Rapoport-Zink spaces to the moduli spaces of global $G$-shtukas. The morphism is finite by the Proposition 6.19.
The author also gives an application of the Main Theorem 0.1 to the leaves inside a Newton stratum and compute dimensions of these leaves which turns out to be the same for all leaves.

\newpage
\end{document}